# Generalized Seikkala Differentiability and its Application to Fuzzy initial value problem

U. M. Pirzada

*School of Engineering and Technology, Navrachana University of Vadodara, India, salmap@nuv.ac.in*

**Abstract**

This paper proposes a new generalized Seikkala derivative (gS-derivative) of a fuzzy-valued function. We see that, there are many elementary fuzzy-valued functions which occur frequently as solution of fuzzy differential equation, are not Seikkala differentiable but they are generalized Seikkala differentiable. We discuss some of the properties of proposed differentiability. Using gS-differentiability, we find the solution of fuzzy initial value problem.

## 1 Introduction

To study the solution of fuzzy initial and boundary value problems, we need the concept of differentiability of fuzzy-valued function. The fuzzy differential calculus is developed by different authors, like Dubois (1982)[1], Puri and Ralescu (1983)[7], Kaleva (1987)[5] and Seikkala (1987)[8]. With different types of differentiability concept, theory of fuzzy differential equation is extensively studied. Recently, solutions to fuzzy differential equations under generalized Hukuhara derivative using contractive-like mapping principles is studied by Elder et al. (2015)[2].

Seikkala differentiability is simpler and stronger than the other differentiability. Unfortunately, there are many elementary fuzzy-valued functions which frequently occur in solution of fuzzy differential equations, are not Seikkala differentiable. The purpose of the paper is to generalize the Seikkala differentiability. We see that a larger class of elementary fuzzy-valued functions belong to a generalized Seikkala differentiability (gS-differentiability ) concept. We find the fuzzy solution of a fuzzy initial value problem using gS-differentiability of a fuzzy-valued function.

The paper is organized in following manner. Section 2 contains basics of fuzzy numbers. The generalized Seikkala differentiability of a fuzzy-valued function is proposed and compared with Seikkala differentiability with appropriate examples in Section 3. The properties of gS-differentiability are discussed in the same section. Conclusion is given in the following section.

## 2 Fuzzy numbers

Fuzzy numbers have been introduced by Lotfi A. Zadeh to handle nonspecific numerical quantities in a practical way. For instance, profit on some product is approximately known, 2 Rs, say, then we can express this approximate amount 2 by means of fuzzy number. The numbers which are near to 2, can also included in fuzzy number 2 with varying membership grade. That is, 1 is regarded as approximately 2 with degree of membership 0.5. Similarly, 3 can be regarded as approximately 2 with degree of membership 0.6. In general, when numbers are imprecise/inexact/partially known, they can be represented as fuzzy numbers rather than real numbers. The mathematical definition of fuzzy number is given as follows.

**Definition 2.1** Let R be the set of real numbers and a: R $\to$ [0, 1] be a membership function of a fuzzy set. Then a is said to be a fuzzy number if there is some $r_0 \in$ R such that $a(r_0) = 1$; function a is quasi-concave and upper semi-continuous on R and clousure of a set $\{r \in R / a(r) > 0\}$ forms a compact set. The set of fuzzy numbers on R is denoted by F(R).

To study the properties of fuzzy numbers, we define corresponding crisp (ordinary) sets.

**Definition 2.2** The α-level set $a_\alpha$ of any a $\in$ F(R) is a crisp set $a_\alpha = \{r \in R / a(r) \geq \alpha\}$, for $\alpha \in (0,1]$ and the 0-level set $a_0$ is defined as the closure of the set $\{r \in R / a(r) > 0\}$.

Using the theorem of Goetschel and Voxman[4], we have following the characterization of a fuzzy number.

**Theorem 2.1** A fuzzy number a is defined by any pair a = $(a_1(\alpha), a_2(\alpha))$ of functions $a_i$ : [0,1] $\to$ R, defining the end-points of the α-level sets, satisfying the following conditions:

(i) $a_1$ is a bounded nondecreasing left-continuous function for $\alpha \in (0,1]$ and right-continuous at $\alpha = 0$;

(ii) $a_2$ is a bounded nonincreasing left-continuous function for $\alpha \in (0,1]$ and right continuous at $\alpha = 0$;

(iii) $a_1(\alpha) \leq a_2(\alpha)$ for all α.

Addition, subtraction, multiplication of two fuzzy numbers are not the same as defined for real numbers. As fuzzy numbers are defined by membership functions, their arithmetic operations can be defined using continuous operation.

**Definition 2.3** Addition and scalar multiplication using their α-level sets are given as follows:

$(a + b)_\alpha = [(a+b)_1(\alpha), (a+b)_2(\alpha)]$, where $(a+b)_1(\alpha) = a_1(\alpha) + b_1(\alpha)$ and $(a+b)_2(\alpha) = a_2(\alpha) + b_2(\alpha)$

$(\lambda a)_\alpha = [(\lambda a)_1(\alpha), (\lambda a)_2(\alpha)]$, if $\lambda \geq 0$
$\quad\quad = [(\lambda a)_2(\alpha), (\lambda a_1)(\alpha)]$, if $\lambda < 0$, where $(\lambda a)_1(\alpha) = \lambda a_1(\alpha)$ and $(\lambda a)_2(\alpha) = \lambda a_2(\alpha)$
for α-level sets $a_\alpha = [a_1(\alpha), a_2(\alpha)]$, $b_\alpha = [b_1(\alpha), b_2(\alpha)]$, for $\alpha \in [0, 1]$ and scalar $\lambda \in R$.

## 3 Generalized Seikkala differentiability and its properties

Rate of change of function is a derivative of that function. If we consider rate of change of imprecise or fuzzy function, we have fuzzy derivative concept. For example, some object is falling from specific height. The motion of the object after t seconds is known approximately can be expressed by fuzzy function. To find instantaneous rate of change of object with respect to time is given as fuzzy derivative at that time. Many authors have defined fuzzy derivatives mathematically in literature. In this paper, we use concept of Seikkala derivative of fuzzy function to study the fuzzy growth and decay systems.

### 3.1 Seikkala differentiability

Seikkala differentiability of fuzzy-valued function y: I $\rightarrow$F(R) is defined as follows. The definition is adopted from Seikkala[3].

**Definition 3.1** Let I be subset of R and y be a fuzzy-valued function defined on I. The α-level sets $y_\alpha(t) = [y_1(t, \alpha), y_2(t, \alpha)]$ for $\alpha \in [0,1]$ and $t \in I$. We assume that derivatives of $y_i(t, \alpha)$, i = 1, 2 exist for all $t \in I$ and for each α.
We define $(y'(t))_\alpha = [y_1'(t, \alpha), y_2'(t, \alpha)]$ for all $t \in I$, all α.

If, for each fixed $t \in I$, $(y'(t))_\alpha$ defines the α-level set of a fuzzy number, then we say that Seikkala derivative of y(t) exists at t and it is denoted by fuzzy-valued function y'(t).

The Seikkala derivative involves two steps:

(1) Check both level functions are differentiable or not
(2) Check level sets of derivatives define fuzzy numbers or not.

Exponential functions are used to represent real-world applications, such as bacterial growth/decay, population growth/decline. If the initial population is imprecise or inexact then to represent the system following fuzzy functions are used. And to find rate of change of such function, we use fuzzy derivative. In particular, Seikkala derivative. So we consider some illustrations of how to find Seikkala derivatives of functions those are useful in real life situations.

The following fuzzy function is exists in fuzzy decay problem. We study the rate of change of the function using uncertain (fuzzy) derivative.

**Example 3.1** Consider a fuzzy-valued function g(t) = a exp(-t), t $\in$ R and a is a fuzzy number with α-level sets $g_\alpha$(t) = [$g_1$(t, α), $g_2$(t, α)] = [$a_1$(α) exp(-t), $a_2$(α) exp(-t)]. To check Seikkala differentiability of given fuzzy-valued function, first we check both its level functions are differentiable or not.

We see that $g_1$(t, α) = $a_1$(α) exp(-x) and $g_2$(x, α) = $a_2$(α) exp(-t) are differentiable for t $\in$ R. Next, we check that the level sets
$$(g'(t))_\alpha = [g_1'(t, \alpha), g_2'(t, \alpha)] = [-a_1(\alpha)\exp(t), -a_2(\alpha)\exp(t)]$$
define a fuzzy number for each t in R. By checking sufficient conditions for $(g'(t))_\alpha$ to define α-level sets of fuzzy number,
(i) $g_1$'(t, α) is an increasing function of α for each t$\in$ R;
(ii) $g_2$'(t, α) is a decreasing function of α for each t$\in$ R; and
(iii) $g_1$'(t, α) ≤ $g_2$'(t, α) for all t $\in$ R, we see that
$\partial g_1$'(t, α)/$\partial \alpha$ = -$a_1$'(α) exp(t) < 0 as $a_1$'(α) > 0 & $\partial g_2$'(t, α)/$\partial \alpha$ = -$a_2$'(α) exp(t) < 0 as $a_2$'(α) < 0.

Therefore, Seikkala derivative of g does not exist.

We consider another example of derivative of fuzzy function which occur in uncertain periodic motion of an object.

**Example 3.2** Consider a fuzzy-valued function h(t) = a sin(t), t in [0, π], where a is a fuzzy number. The α- level sets of h(t) are [$a_1$(α)sin(t), $a_2$(α)sin(t)]. The level function are differentiable but their derivatives $h_1$'(t, α) = $a_1$(α) cos(t) and $h_2$'(t, α) = $a_2$(α) cos(t) does not define fuzzy number for each t in [π/2, π] and hence h is not Seikkala differentiable for t in [π/2, π].

**3.2 Generalized Seikkala differentiability**

From Example 3.1 and 3.2, we see that not all fuzzy-valued functions are Seikkala differentiable. The generalized Seikkala derivative (gS-derivative) of a fuzzy-valued function is defined as follows:

**Definition 3.2** Let I be a real interval. A fuzzy-valued function f: I → F(R) with α-level sets $f_\alpha$(t) = [ $f_1$(t, α), $f_2$(t, α)], for t $\in$ I and α $\in$ [0, 1] is said to have generalized Seikkala derivative f '(t) if $f_1$(t, α) and $f_2$(t, α) are differentiable for each t $\in$ I and
f '$_\alpha$(t) = [ min{$f_1$′(t, α), $f_2$′(t, α)}, max{ $f_1$′(t, α), $f_2$′(t, α)}],
for all α defines a fuzzy number for each t in I.

We see that the uncertain functions defined in Example 3.1 and 3.2 exist in real life, are gS-differentiable.

**Example 3.3** The fuzzy-valued function g(t) = a exp(-t), t in R, defined in Example 3.2. The derivatives of level functions of g(t) are $g_1$'(t, α) = -$a_1$(α) exp(-t) and $g_2$'(t, α) = -$a_2$(α) exp(-t). By definition of gS-differentiability, α-level sets $g_\alpha$'(t) defined as

$g_\alpha$'(t) = [min{-$a_1$(α) exp(-t), -$a_2$(α) exp(-t)}, max{-$a_1$(α) exp(-t), -$a_2$(α) exp(-t)}]

which is equal to $g_\alpha$'(t) = [-$a_2$(α) exp(-t), -$a_1$(α) exp(-t)]

as - $a_2(\alpha) \leq$ - $a_1(\alpha)$ and exp(-t) $\geq$ 0 for all t.

Therefore, g is gS-differentiable with derivative g′(t) = -a exp(-x).

**Example 3.4** The fuzzy-valued function h(t) defined in Example 3.3 is gS-differentiable with derivative h′(t) = a cos(t). The α-level sets of h′(t) are
[$a_1(\alpha)$ cos(t), $a_2(\alpha)$ cos(t)] for t ∈ [0, π/2] and [$a_2(\alpha)$ cos(t), $a_1(\alpha)$ cos(t)] for x ∈ (π/2, π].

### 3.3 Properties

We prove some properties of gS-differentiability.

**Theorem 3.1** If fuzzy-valued function f: I → F(R) is S-differentiable then it gS-differentiable.

**Proof.** Let f be a S-differentiable function then $f_\alpha$'(t) = [$f_1$′(t, α), $f_2$′(t, α)] defines a fuzzy number for each t in I. We can write f ′$_\alpha$(t) = [$f_1$′(t, α) = min{$f_1$′(t, α), $f_2$′(t, α)}, $f_2$′(t, α) = max{ $f_1$′(t, α), $f_2$′(t, α)}]
which defines fuzzy number for each t in I. Hence the proof.

**Remark 3.1** Converse of the above theorem is not true. That is, if f is gS-differentiable it may not be S-differentiable. For instance, functions defined in Example 3.3 and 3.4 are gS-differentiable but they are not S-differentiable.

Another property is regarding relation between gS-differentiability of fuzzy-valued function and differentiability of α-level functions.

**Theorem 3.2** The fuzzy-valued function f: I → F(R) is gS-differentiable if and only if its α-level functions $f_1$(t, α) and $f_2$(t, α) are differentiable for t in I and for each α.

**Proof.** The result followed by definition of gS-differentiability of fuzzy-valued function.

**Theorem 3.3** If fuzzy-valued functions f and g are gS-differentiable then f+g is gS-differentiable.

**Proof.** Let f and g be gS-differentiable function then
f ′$_\alpha$(t) = [min{$f_1$′(t, α), $f_2$′(t, α)}, max{ $f_1$′(t, α), $f_2$′(t, α)}] and
g ′$_\alpha$(t) = [min{$g_1$′(t, α), $g_2$′(t, α)}, max{ $g_1$′(t, α), $g_2$′(t, α)}]

(f+g) is gS-differentiable if (f+g)'$_\alpha$(t) given as
(f+g)'$_\alpha$(t) = [min{(f+g)$_1$′(t, α), (f+g)$_2$′(t, α)}, max{ (f+g)$_1$′(t, α), (f+g)$_2$′(t, α)}] defines a fuzzy number for each t.

Here, (f+g)$_1$'(t, α) = $f_1$' + $g_1$' and (f+g)$_2$'(t, α) = $f_2$' + $g_2$' if $f_1$' ≤ $f_2$' and $g_1$' ≤ $g_2$'
    (f+g)$_1$'(t, α) = $f_1$' + $g_2$' and (f+g)$_2$'(t, α) = $f_2$' + $g_1$' if $f_1$' ≤ $f_2$' and $g_2$' ≤ $g_1$'
    (f+g)$_1$'(t, α) = $f_2$' + $g_1$' and (f+g)$_2$'(t, α) = $f_1$' + $g_2$' if $f_2$' ≤ $f_1$' and $g_1$' ≤ $g_2$'
    (f+g)$_1$'(t, α) = $f_2$' + $g_2$' and (f+g)$_2$'(t, α) = $f_1$' + $g_1$' if $f_2$' ≤ $f_1$' and $g_2$' ≤ $g_1$'.

In all these four cases, f+g is gS-differentiable as f and g are gS-differentiable.

**Theorem 3.4** If fuzzy-valued function f is gS-differentiable and $\lambda$ in R then ($\lambda$f) is gS-differentiable.

**Proof.** Similar arguments as in Theorem 3.3, the result can be proved.

## 4 Application

Now we consider an application of gS-derivative to solve a fuzzy initial value problem. The first example is real life example of decay of radioactive isotope where the initial amount is uncertain. Thus the system is studied as fuzzy system because uncertainty exists in the system. The mathematical modeling of such problem gives us fuzzy initial value problem. To find rate of change of approximate amount of radioactive isotope at any time t, we use gS-differentiability concept.

**Example 4.1** A radioactive isotope has initial mass is around 200 mg, which is decreasing at a rate of present mass. Find the expression for the approximate amount of isotope remaining at any time.

**Solution: Fuzzy Model**

In this case, we have a fuzzy initial value problem

$$dy/dt = f(t, y) = -y \quad \text{(rate of decay is equal to amount present)} \text{------ (1)}$$

$y(0)$ = around 200, $x \in I = [0, 1]$, y is a fuzzy-valued function on I, f is a fuzzy-valued function on $I \times F(R)$ and c = (around 200) is a fuzzy number with $\alpha$-level sets $c_\alpha = [c_1(\alpha), c_2(\alpha)]$.

**Solution**

If we assume that y is S-differentiable with S-derivatives. The equation (1) can be rewritten as

$$y_1'(t, \alpha) = -y_2(t, \alpha) \text{-------- (2)}$$
$$y_2'(t, \alpha) = -y_1(t, \alpha) \text{-------- (3)}$$

with initial conditions $y_1(0, \alpha) = c_1(\alpha)$ and $y_2(0, \alpha) = c_2(\alpha)$, for $\alpha \in [0, 1]$.

Solving equations (2) and (3), we get

$y_1(t, \alpha) = 0.5(c_1(\alpha) - c_2(\alpha))\exp(t) + 0.5(c_1(\alpha) + c_2(\alpha))\exp(-t)$, ----- (4) and
$y_2(t, \alpha) = 0.5(c_1(\alpha) + c_2(\alpha))\exp(-t) - 0.5(c_1(\alpha) - c_2(\alpha))\exp(t)$, ----- (5) for $\alpha \in [0, 1]$.

If y is Seikkala differentiable and $[y_1, y_2]$ defines fuzzy number, for each $t \in I$, then y will be fuzzy solution of the problem. But as we have discussed in Pirzada(2017)[6], y is not Seikkala differentiable and therefore, it cannot be a solution of problem (1).

Now we consider the given problem under gS-differentiability, we see that y(t) is gS-differentiable with

$y'_\alpha(x) = [y_2'(t, \alpha), y_1'(t, \alpha)]$, as sufficient conditions for existence of a fuzzy number y' for each x:

$\partial y_2'(t, \alpha) / \partial \alpha > 0$
and
$\partial y_1'(t, \alpha) / \partial \alpha < 0$,

for given fuzzy initial condition y(0) = c = (around 200) = (195, 200, 205). Hence equations (4) and (5) define fuzzy solution of the problem.

**Interpretation**

The amount of isotope at any time t is represent approximately using functions $y_1(t, \alpha)$ and $y_2(t, \alpha)$, $\alpha \in [0,1]$ given in equations (4) and (5) respectively. That is, we have infinite family of solution $(y_1, y_2)$. For $\alpha = 0.5$, solution can be represented by two curves given following figure.

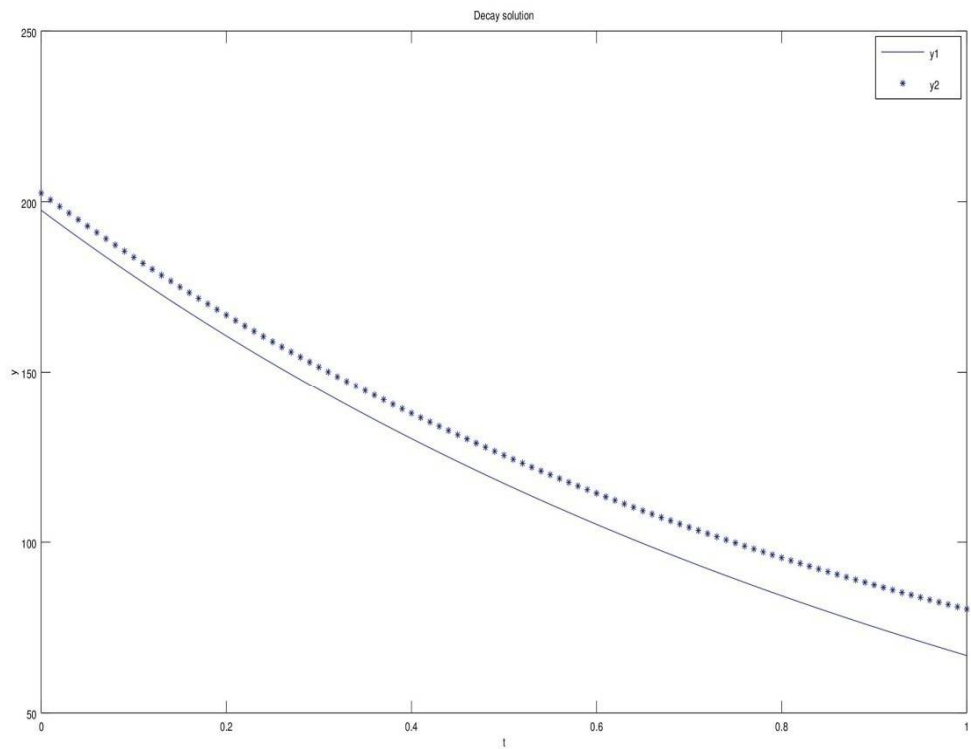

**Example 4.2** Suppose biochemical oxygen demand (BOD) in water is of amount approximately $D_0$. The rate of decay is proportional to the amount of dissolved oxygen in water at present. Determine the approximate amount of dissolved oxygen in water at time t.

**Solution:** Let the approximate amount of dissolved oxygen is D(t), D(t) is a fuzzy-valued function of time t. The problem governed by the above situation is
$$dD / dt = - D, D(0) = D_0,$$
where $D_0$ is a fuzzy number (90, 100, 120) and negative sign indicate decay. Assume that D is gS-differentiable and therefore the fuzzy equation can be rewiritten as

$$dD_1(t,\alpha) / dt = -D_2(t,\alpha), \quad dD_2(t,\alpha) / dt = -D_1(t,\alpha)$$

with initial condition: $D_1(0, \alpha) = D_{01}(\alpha)$ and $D_2(0,\alpha) = D_{02}(\alpha)$. Solving these equations, we get

$D_1(t, \alpha) = 0.5(c_1(\alpha) - c_2(\alpha))\exp(t) + 0.5(c_1(\alpha) + c_2(\alpha))\exp(-t)$, and
$D_2(t, \alpha) = 0.5(c_1(\alpha) + c_2(\alpha))\exp(-t) - 0.5(c_1(\alpha) - c_2(\alpha))\exp(t)$, for $\alpha \in [0, 1]$. The approximate amount of dissolved oxygen in water at time t is a family of solutions $(D_1, D_2)$.

**Conclusions**

The generalization of Seikkala derivative is defined in the proposed work. Different properties are discussed. We observed that the gS-differentiability is useful in finding the solution of fuzzy initial value problem. We have seen that solution of fuzzy decay model exists under gS-differentiability only. The solution of the problem does not exist under Seikkala differentiability of fuzzy-valued function. Further, we can extend its application to solve fuzzy boundary value problems and fuzzy partial differential equations.